\documentclass[l10 pt]{article}

\usepackage{graphicx,graphics,epsfig,subfigure}
\usepackage{times,color}
\usepackage{amssymb,amsfonts,amsmath,latexsym,amscd}
\usepackage{epstopdf}
\newcommand{\grad}{\nabla}
\usepackage{graphics,epsfig}

\newcommand{\eps}{\varepsilon}
\newcommand{\kerc}{c^\bot}
\newcommand{\cone}{{a}}

\title{\LARGE \bf Root locii for systems defined on Hilbert spaces}

\author{
Birgit Jacob \thanks{B. Jacob is with the Faculty of  Mathematics und Natural Sciences, University of Wuppertal,  Gau\ss stra\ss e 20, 42119 Wuppertal, GERMANY,  Email: {\tt\small bjacob@uni-wuppertal.de} . }
\and
Kirsten Morris\thanks{K. Morris is with the University of Waterloo,
Department of Applied Mathematics, Waterloo, ON, N2L 3G1,
CANADA, Email: {\tt\small kmorris@uwaterloo.ca}.}
\thanks{The authors gratefully acknowledges financial support for this research
from the National Sciences and Engineering Research Council, Deutsche Forschungsgemeinschaft (Grant JA 735/9-1) and the Banff International Research Station.
}
}

\usepackage{amssymb,amsmath}

\newtheorem{thm}{Theorem}[section]
\newtheorem{lem}[thm]{Lemma}

\newtheorem{prop}[thm]{Proposition}

\newtheorem{cor}[thm]{Corollary}
\newtheorem{defn}[thm]{Definition}
\newtheorem{eg}[thm]{Example}

\newtheorem{conj}[thm]{Conjecture}

\newenvironment{mat}{ \left[ \begin{array}{ccccc} }{\end{array} \right]  }
\newcommand{\mc}{{\mathbb C}}
\newcommand{\re}{\mathcal{R}}

\newcommand{\cL}{\cal L}
\newcommand{\cZ}{{\cal Z}}

\newcommand{\cl}{{\mathcal L}}

\newcommand{\dom}{{\cal D}}

\newcommand{\la}{\langle}
\newcommand{\ra}{\rangle}

\newcommand{\half}{\frac{1}{2}}
\newcommand{\inv}{{\rm inv}}

\newcommand{\be}{\begin{equation}}
\newcommand{\ee}{\end{equation}}
\newcommand{\bea}{\begin{eqnarray}}
\newcommand{\eea}{\end{eqnarray}}
\newcommand{\beann}{\begin{eqnarray*}}
\newcommand{\eeann}{\end{eqnarray*}}

\renewcommand{\re}{{\rm Re}\,}

\newcommand{\pf}{ \noindent {\em Proof:}  }
\newcommand{\eop}{{\hspace*{\fill}$\square$}\vspace{1ex}}

\renewcommand{\conj}[1]{\overline{#1}}
\renewcommand{\bar}[1]{\overline{#1}}
\newcommand{\Di}{D^{-1}}

\begin{document}
\maketitle
\begin{abstract}
The root locus is an important tool for analysing the stability and time constants of linear finite-dimensional systems as a parameter, often the gain, is varied. However, many systems are modelled by partial differential equations or delay equations. These systems evolve on an infinite-dimensional space and their transfer functions are not rational. In this paper a rigorous definition of the root locus for infinite-dimensional systems is given and it is shown that the root locus is well-defined  for a large class of infinite-dimensional systems. As for finite-dimensional systems, any  limit point of a branch of the root locus is a zero. However, the asymptotic behaviour can be quite different from that for finite-dimensional systems. This point is illustrated with a number of examples. It is shown that the familiar pole-zero interlacing property for collocated systems with a Hermitian state matrix extends to infinite-dimensional systems with self-adjoint generator. This interlacing property is also shown to hold for collocated systems with a skew-adjoint generator. 
\end{abstract}
\section{Introduction}

Consider the  control system on a Hilbert space $\cZ$ 
\begin{equation*}
\begin{array}{lll}
\dot{x}(t) &=& Ax(t)+Bu(t)\\
y(t) &=& Cx(t).
\end{array}
\end{equation*}
where for some $b$, $c \in \cZ$, $B u = b u$, $C x = \la c ,x \ra $.
The simplest control for a system is a constant gain,
 $$u(t) =-k y(t) + v(t),$$ where  $k>0$ is real and $v(t)$ is an external signal.
 Thus,   the eigenvalues of  $A-k BC$ as $k\rightarrow \infty$ are  of interest.  A plot of these eigenvalues as $k \rightarrow \infty$ is known as a root locus plot.
An understanding of the behaviour of these eigenvalues as $k$ varies, or the root locus, is important to understanding the behaviour of the system with feedback.

Suppose that the system is finite-dimensional; that is $A\in \mathbb C^{n\times n}$.  If the relative degree of the system is $r$  then there are $r$ eigenvalues going to infinity and the remaining eigenvalues tend to the zeros of the transfer function 
 \cite[e.g.]{MacFarlane77,KouvaritakisMacFarlane76}. 
 Furthermore, the angle of the asymptotes as $k \rightarrow \infty$, in particular whether they are in the left-half-plane, is determined  by the relative degree.  

Extension of these results, now well-known for finite-dimensional systems, to in\-finite-dimensions has been elusive. In  \cite{Pohjolainen1981} the root locus is considered for the case where  $A$ is self-adjoint with compact resolvent on a Hilbert space $\cZ$, $B$ is a linear bounded operator from $\mathbb C^p$ to $\cZ$, and $C: D(C) \rightarrow \mathbb C^p$ where $D(A) \subset D(C)$ is $A$-bounded.  A complete analysis of collocated boundary control of parabolic systems on an interval was provided in \cite{ByrnesGilliam}. The analysis in that paper uses  results from differential equations theory and is difficult to extend to  more general classes.
In \cite{Kobayashi92} high-gain  output feedback of infinite-dimensional systems in the case where $A$ generates an analytic semigroup  and $B=C^*$ was studied. 
The zeros of the system are given as the eigenvalues of an operator and a nonlinear stabilizing feedback law is constructed. Zeros of systems where $A$ is self-adjoint and $B=C^*$  are shown to be real and  be bounded by $\alpha$ if $A+A^* \leq 2 \alpha I $ on $D(A)$ in \cite{ZwartHof}. If moreover, the system transfer function can be written in spectral form, and additional technical conditions are satisfied, the poles and the zeros interlace on the real axis.


A review of the definitions and properties of the zeros for systems with bounded control and observation is first presented.  It is shown that in many situations, the zeros are in the  pseudo-spectrum of the generator. The application of the pseudo-spectrum to the analysis of zeros is new. This is used to show that  in many cases,   even if no zeros coincide with an eigenvalue,   the zeros become asymptotically close to the eigenvalues. This extends an earlier result \cite[Thm. 4.4]{ZwartHof} that was obtaining using a different approach.
In this paper the root locus for single-input-single-output  infinite-dimensional systems is shown to be well-defined. If no invariant zeros are in the spectrum of $A$ each eigenvalue of $A$ defines a  branch of the root locus and these curves are smooth and non-intersecting.  Moreover, if any branch converges to a point, that point is a zero of the system. Conversely, each zero is the terminus of a branch of the root locus.  The root locus of systems where $B=C^*$ and the generator $A$ is either self-adjoint or skew-adjoint is considered in detail. It is shown that the zeros interlace with the poles, that the root locus and the root locus is contained in the left half-plane. Although both self-adjoint and skew-adjoint systems are relative degree one when $B=C^*$, the asymptotic behavior of the root locus is very different.
Preliminary  versions of  some of the results in sections III and V appeared in a CDC paper, reference \cite{JacoMorrCDC}. 

\section{Zeros}

 Consider control systems $\Sigma(A,B,C,D)$ on a Hilbert space ${\cal Z}$ of the form
\begin{equation}
\begin{array}{lll}
\dot{x}(t) &=& Ax(t)+Bu(t)\\
y(t) &=& Cx(t)+Du(t),
\end{array}
\label{sys}
\end{equation}
where $A$ is the generator of $C_0$-semigroup on  $\cZ$, $B$ and $C$ are bounded operators with scalar input and output spaces, that is, $Bu=b u$ and $Cx = \langle x,c  \rangle$ for some $b,c \in \cZ$, $D\in \mathbb C$, and  $B, C\not= 0$. We write $\Sigma(A,B,C):= \Sigma(A,B,C,0)$ for short.

 We denote by $D(A)$, $\sigma(A)$, and $\rho(A)$,
the domain, the spectrum, and the resolvent set, respectively. It is assumed throughout this paper that  $\sigma(A)$ is non empty and consists of isolated eigenvalues with finite algebraic multiplicity only. That is $ \sigma(A)$ is finite or countable, with no finite accumulation point, consisting of eigenvalues of $A$ of  finite algebraic multiplicity. Equivalently,   
$\sigma_{{\rm ess}}(A)\not=\emptyset$. Here $\sigma_{{\rm ess}}(A)$ denotes  the set of all $\lambda \in \mathbb C$ such that $\lambda I-A$ is not a Fredholm operator and an operator $T$ on $\cZ$ is called a {\em Fredholm operator} if $T$ is closed and the dimension of the kernel of  $T$ and the dimension of $\cZ / {\rm Im} T$ are finite. 
The assumption on $\sigma(A)$ is for example satisfied if $A$ has compact resolvent.

Let $\{\lambda_n\}$ be the set of eigenvalues of $A$  counted according to their multiplicity.

\begin{defn}
Let $(T(t))_{t\ge 0}$ be the $C_0$-semigroup generated by $A$.
Then the system $\Sigma(A,B,C,D)$ is called {\em approximately observable} if for every $z\in \cZ\backslash\{0\}$ the function $CT(t)z$ is not identically zero on $[0,\infty)$.

The system $\Sigma(A,B,C,D)$ is called {\em approximately controllable} if the system $\Sigma(A^*, C^*, B^*,D^*)$ is approximately observable.
\end{defn}

Let $G(s) =C (sI-A)^{-1}B+D$, $s\in \rho(A)$, indicate the transfer function of the system $\Sigma(A,B,C,D)$, defined using the characteristic function \cite{Zwart2004}.
If $\Sigma(A,B,C,D)$ is approximately observable, this definition is equivalent to the other definitions of the transfer function for all $s \in \rho (A)$ \cite[Cor. 2.8]{Zwart2004}. 

Transmission zeros and invariant zeros are defined similarly to the finite-dimensional case. 
 
\begin{defn}
A complex number $s\in \rho(A)$ is a  {\em transmission zeros of $\Sigma(A,B,C,D)$} if
  $G(s) = C (sI-A)^{-1} B+D=0$. 
\end{defn}

\begin{defn}
The {\em invariant zeros of $\Sigma(A,B,C,D)$} 
are the set of all $\lambda\in\mathbb C$ such that
  \begin{equation}
      \label{star}
\begin{mat} \lambda I - A & -B \\ C & D \end{mat} \begin{mat} x \\ u
\end{mat}
= \begin{mat} 0 \\ 0 \end{mat}
\end{equation}
has a solution for  some scalar $u$ and non-zero $x \in D(A)$.
Denote the set of invariant zeros of a system by $\inv (A,B,C,D)$.
\end{defn}

As in finite-dimensional systems, it is straightforward to show
that   $\lambda \in  \rho(A)$ is a transmission zero if and only it is an invariant zero, see \cite{ZwartHof} for the case $D=0$.
 
\begin{defn}
The system $\Sigma(A,B,C,D)$ is called {\em minimal} if $\inv (A,B,C;D)\cap \sigma(A) =\emptyset$.
\end{defn}

Minimality will be important throughout this paper. 
The term {\em minimal} is often used in finite-dimensional systems theory to refer to a system that is controllable and observable. 
An operator $A$ is {\em Riesz-spectral} if $A$ is closed, has only simple eigenvalues $\{\lambda_n\}_n$, the corresponding eigenvalues form a Riesz basis of $\cZ$, and the closure of the point spectrum is totally disconnected.

\begin{prop}\label{conj1}
If $A$ is Riesz-spectral and  $\Sigma(A,B,C,D)$ is approximately  observable and  approximately controllable   then the system is minimal. 
 \end{prop}

\pf  See \cite[Ex. 4.28b]{CZbook} for $D=0$. The proof is similar to the proof for  finite-dimensional systems and the generalization to $D\not=0$ is straightforward. \eop

\begin{lem}
\label{lem:Gk}
 \cite[Theorem 1.2]{WeissXu} Let  $\Sigma(A,B,C,D)$ be a system. We assume $k\not=0$ and  $kD\not=-1$. Then for
 any point $s \in \rho(A)$ and $k \neq 0$, we have  $s \in \sigma (A-B k (1+Dk )^{-1} C)$ if and only if $G(s) =-\frac{1}{k}$.
\end{lem}

\begin{thm}
\label{thm:allk}
Suppose that the system $\Sigma(A,B,C,D)$ with $kD\not=-1$ is either minimal or  approximately controllable and approximately observable,  then 
\[\sigma (A) \cap \sigma (A-B k (1+Dk )^{-1} C)=\emptyset \] for all $k\not=0$.
Moreover, we have $\sigma (A-B k_1 (1+Dk_1 )^{-1} C)\cap \sigma (A-B k_2 (1+Dk_2 )^{-1} C)=\emptyset$ for $k_1\not=k_2$.
\end{thm}

\pf
The proof  for approximately controllable and approximately observable systems can be found  in \cite[Theorem 1.4]{WeissXu}. Assume now that the system is minimal. 
Suppose that for some $k \neq 0$ with $kD\not=-1$ there is $s \in  \sigma(A-B k (1+Dk )^{-1} C) \cap \sigma (A) $ and let $x_o \neq 0 $, $x_k \neq 0$ be such that
\begin{equation} s x_o = A x_o,  \hspace{2em} s x_k = A x_k - B k (1+Dk )^{-1} C x_k . \label{eqn27} \end{equation}
Assume $Cx_k=0$. Then $s$ is an invariant zero of $\Sigma(A,B,C,D)$, which is a contradiction to the minimality of the system.
Thus there exists a scalar  $\alpha$  such that 
$$ C x_o + \alpha ( 1- Dk(1+Dk)^{-1}) C x_k =0 . $$
If $x_o + \alpha x_k =0$, then  from (\ref{eqn27}) it follows that $C x_{k}  =0$. Here we also used the fact that $B\not =0$. 
Suppose now that $x_o + \alpha x_k \neq 0$.
Then  
$$
\begin{mat} s I - A & -B  \\ C & D \end{mat} \begin{mat} x_o +  \alpha x_k \\ -\alpha k(1+Dk)^{-1})  C x_k
\end{mat}
= \begin{mat} 0 \\ 0 \end{mat} 
$$
and so $s$ is an invariant zero of $\Sigma(A,B,C,D)$. This implies again that  $\sigma (A) \cap \inv (A,B,C,D) $ is not empty. 
Finally, $\sigma (A-B k_1 (1+Dk_1 )^{-1}C)\cap \sigma (A- Bk_2 (1+Dk_2 )^{-1} C)=\emptyset$ for $k_1\not=k_2$ follows from Lemma \ref{lem:Gk}.
\eop

\begin{prop}\label{prop:mero}
Suppose that the system $\Sigma(A,B,C,D)$ is minimal. Then the transfer function $G(s)$ of $\Sigma(A,B,C,D)$ is meromorphic on $\mathbb C\backslash \sigma(A)$ and each $\lambda \in \sigma(A)$ is a pole of $G(s)$.
\end{prop}
\pf  Because $\sigma(A)$ consists of isolated eigenvalues with finite algebraic multiplicity,   $(sI-A)^{-1}$ is meromorphic on  $\mathbb C$ and each $\lambda \in \sigma(A)$ is a pole \cite[III.6.5, pg. 180]{Kato}. Thus $G(s)$ is meromorphic as well and each $\lambda \in \sigma(A)$ is a pole or a removable singularity. It remains to show that each $\lambda \in \sigma(A)$ is a pole of $G(s)$. 
Let $\lambda \in \sigma(A)$. Then there is a sequence  $\{s_k\}_{k\in \mathbb N}$ converging to $\lambda$ such that each $s_k$ is an eigenvalue of $A-\frac{1}{k}B  (1+D\frac{1}{k} )^{-1}C$, see \cite[IV.3.5, pg.~213]{Kato}. Theorem \ref{thm:allk} implies that $s_k\in \rho(A)$ for every $k\in \mathbb N$ and by Lemma \ref{lem:Gk} we obtain $G(s_k) =-k$. Thus $\lambda$ is a non-removable singularity of $G(s)$.\eop

In particular, the previous proposition shows that since  $\sigma(A)\not=\emptyset$, the transfer functions of minimal systems are not identically zero.

\begin{prop}
\label{conj_invzeros}
Let $\Sigma(A,B,C,D)$ be a minimal system.  Then the set $\inv (A,B,C,D)$ is countable and has no finite accumulation point.
\end{prop}

\pf
Since $G$ is not identically zero and meromorphic, the set of transmission zeros of $\Sigma(A,B,C,D)$ is countable and  has no  finite accumulation point. Moreover, by assumption the spectrum of $A$  is countable and  has no finite accumulation point. 
Since the set $\inv (A,B,C,D)$ is a subset of the union of  $\sigma(A)$ and the set of transmission zeros of $\Sigma(A,B,C,D)$, $\inv (A,B,C,D)$ is countable and  has no  finite accumulation point. 
\eop

The invariant zeros can  be characterized as the spectrum of an operator. The case where $D \neq 0$ is simplest and is presented first. 

\begin{prop}
The invariant zeros of $\Sigma(A,B,C,D)$  with $D\not=0$ are the eigenvalues of the operator $A_\infty=A+B\Di C$ with $D(A_\infty)=D(A)$.
\end{prop}

\pf If  $\lambda\in \mathbb C$ is an invariant zero of $\Sigma(A,B,C,D)$, that is, there exists $x\in D(A)\backslash\{0\}$ and $u\in\mathbb C$ such that  $\lambda x-Ax-Bu=0$ and $Cx+Du=0$, then $u=-D^{-1}Cx$ and therefore $(A-B\Di C)x=\lambda x$. Thus $\lambda$ is an eigenvalue of $A-B\Di C$. The converse implication can be proved in a similar manner.\eop

Indicate the kernel of $C$ by
$$\kerc := \{x\in X \mid \langle
x,c\rangle = 0\}.$$  

\begin{thm} \label{bcne0}
\cite[Thm. 2.3]{MorrisRebarberMCSS}, \cite[Thm. 2.3]{MorrisRebarberIJC}
 Suppose that $\Sigma(A,B,C)$ is a minimal system with $\langle b,c\rangle \neq 0$. For $z \in D(A_\infty) = D(K)= D(A), $ define
\begin{equation}
    \label{Ainfinity} 
    Kz=-\frac{\langle Az, c\rangle}{\langle b,c\rangle}, \quad  A_\infty z= Az + b Kz . 
\end{equation} 
Then
 $(A+bK)(\kerc \cap D(A)) \subset \kerc$  and the invariant zeros of $\Sigma(A,B,C)$ are  eigenvalues of $A_{\infty} \vert_{\kerc}$. Moreover, denoting by $\{\mu_n\}$ the invariant zeros of $\Sigma(A,B,C)$, the corresponding eigenfunctions of  $A_{\infty} \vert_{\kerc}$ are given by $\{ (\mu_n I-A)^{-1}b\}$.
\end{thm}

The operator $A_\infty$ is  unique up to addition  of another operator $ b K$ where $Kz=0$ for   $ z\in \kerc $.
If $ c \in D(A^*)$, the perturbation of $A$ is bounded, but in general $A_\infty$ involves an unbounded perturbation and may not generate a $C_0$-semigroup  \cite{MorrisRebarberMCSS}.

If $c \in D(A^{*n})$ for
some integer $n \geq 1$, define
$$Z_{n} = \kerc \cap ( A^{*}c)^{\perp} \cap \cdots
(A^{*n}c)^{\perp},$$
 and define $Z_0 = \kerc$ and $Z_{-1} = X$.  

\begin{thm}\label{bnotin}
\cite[Thm. 2.3]{MorrisRebarberMCSS}, \cite[Thm. 2.7]{MorrisRebarberIJC}
      \label{nested}
      \label{thminv}
 Suppose that $\Sigma(A,B,C)$ is a minimal system and that  an integer $n \geq 1$ exists  such that
\begin{equation}
    \label{bin} c \in D(A^{*n}), ~~~~~~ b \in Z_{n-1}
    \end{equation}
and
\begin{equation}
 \langle b, A^{*n}c\rangle \neq 0.
\end{equation}
The invariant zeros of $\Sigma(A,B,C)$ are  eigenvalues of $A_\infty \vert_{Z_n^\perp} $ where $A_\infty = A+ b K $ and 
\begin{equation}
\label{Kn} Kx =\langle Ax, \cone  \rangle , \hspace{2em} \cone =
\frac{-A^{n*}c}{\langle b, A^{n*}c\rangle}, \hspace{2em} D(K) = D(A).
\end{equation}
\end{thm}
As in the case where $\la b,c \ra \neq 0 $,  changing $K$ on  $(Z_n)^\perp$ does not change the conclusion of Theorem \ref{bnotin}.

The invariant zeros can also be characterized as the eigenvalues of an operator on an invariant subspace in the general case, where $c \notin D(A^{*n}),$ but the definitions are not straightforward because a largest invariant subspace  might not exist. For details, see \cite{MorrisRebarberIJC}.

\begin{defn}
For any $\eps >0$ the {\em $\eps$-pseudospectrum} of an operator $A: D(A)\subset \cZ \rightarrow \cZ$ is
\begin{align*}
&\sigma_\eps (A) \\
&=\{ s \in \mc \mid \| s z - A z \| < \eps\; \mathrm{ for \; some } \; z \in D(A), \|z\|=1 \} .
\end{align*}
\end{defn}

See \cite[pg. 31]{TEbook} for the definition and further properties of the $\eps$-pseudospectrum. In general, the $\eps$-pseudospectrum of an operator  can be quite different from its spectrum. However, for normal operators the $\eps$-pseudospectrum equals the union of $\eps$-balls around the spectrum of $A$.

\begin{thm} \cite{TEbook}
If $A$ is normal, 
$$\sigma_\eps (A) = \bigcup_n B(\lambda_n , \eps ), $$
where $B( \lambda, \eps ):=\{ s\in \mathbb C\mid |s-\lambda|<\eps\}$.
\end{thm}

The following theorem shows that under certain assumptions, the zeros are asymptotically close to the $\eps$-pseudospectrum of $A.$
A sequence $\{\phi_n\}$ in $\cZ$ is called a {\em Riesz system} in $\cZ$ if there exists an isomorphism $S\in L(\cZ)$ such that $\{S\phi_n\}$ is an orthonormal system in $\cZ$.

\begin{thm}
\label{zeros_close}
Suppose that $\Sigma(A,B,C)$ is a minimal system with $\la b,c \ra  \neq 0$, the eigenfunctions of $A_\infty$, see (\ref{Ainfinity}), corresponding to the invariant zeros of $\Sigma(A,B,C)$ form a Riesz system, and $c \in D(A^*) .$  Write the invariant zeros of $\Sigma(A,B,C)$ as $\{ \mu_1 , \mu_2, \ldots\}$ (repeated according to multiplicity) and indicate the corresponding normalized eigenfunctions of $A_\infty$ by $\{ z_1, z_2 , \ldots \}$.  
Then for any $\eps>0$ there is $N$ so that  for all $n>N$   
$$\|A z_n - \mu_n z_n \| < \eps, $$
that is, $ \mu_n\in \sigma_\eps(A)$.
\end{thm}

\pf As $z_n$ is an eigenvector of $A_\infty$ with respect to the eigenvalue $\mu_n$,
\begin{eqnarray*}
\| A z_n - \mu_n z_n \| &\leq& \frac{\| b \|}{\la b,c\ra} | \la A z_n , c \ra | \\
&\leq& \frac{\| b \|}{\la b,c\ra} | \la z_n , A^*  c \ra | .
\end{eqnarray*}
Since $\{z_n\}$ is a Riesz system, they are weakly convergent to $0$ and the result follows.
\eop

If $A$ generates a bounded $C_0$-semigroup, then the fractional powers $(-A)^\alpha$ are well-defined, see \cite{Bala}. 
The following lemma will be useful for establishing several results.

 \begin{lem}
 \label{lem:alpha}
 Suppose $\Sigma(A,B,C)$ is a minimal system,  $A$ generates a bounded $C_0$-semigroup and $\alpha\in \mathbb R$, and  $C (-A)^\alpha $ defines a bounded operator.
 Then the invariant zeros of the systems $\Sigma(A,B,C)$ and $\Sigma(A, (-A)^{-\alpha} B , C (-A)^\alpha  ) $ are identical; and similarly the transmission zeros are identical. 
 \end{lem}
 
 \pf
 If $\mu$ is an invariant zero of $\Sigma(A,B,C)$, then there is $z \in \cZ\backslash \{0\}$ so
$$ \begin{mat} \mu I - A & -B \\ C & 0 \end{mat} \begin{mat} z \\ 1
\end{mat}
= \begin{mat} 0 \\ 0 \end{mat} 
$$
and in fact
$$ z= (\mu I-A)^{-1} B 1 .$$
Defining $\tilde{z}=(\mu I-A)^{-1} (-A)^{-\alpha } B 1 $ yields the result since powers of $-A$ commute with its resolvent.  The equality of the two sets of transmission zeros follows identically.
 \eop

 The previous result can be generalized to the case where $\Sigma(A, (-A)^{-\alpha} B , C (-A)^\alpha  ) $ is a regular system,  using an appropriate generalization of the definition of an invariant zero, but this generalization is not needed in this paper.

\section{Root locus }

In this section it is shown that the root locus of a large class  of infinite-dimensional systems   (\ref{sys}) consists of well-defined curves.  Consider first the case where $D=0. $

The root locus is basically $\sigma (A-k BC)$ for real $k$,  as $k$ moves from $0$ to infinity. 
Thanks to our assumptions on the spectrum of $A$, that is, $\sigma(A)$ is non empty and consists of isolated eigenvalues with finite algebraic multiplicity only,  for every $k\ge 0$ the set $ \sigma(A-kBC)$ consists entirely of isolated eigenvalues  with finite algebraic multiplicity,
 \cite[ XVII Corollary 4.4]{GGK}.
Recall that  the eigenvalues of $A$ are indicated by $\{\lambda_n\}.$ 

There is a family of curves $f_n (k)$ associated to the eigenvalues of $A$ with $f_n (0)= \lambda_n$. The values of $f_n (k)$ are the eigenvalues of  $A-k BC$.   The {\em root locus} is the set of curves $f_n (k)$. In general, it may occur that he root locus is empty, see Example \ref{emptylocus}. However, our assumptions on $\sigma(A)$ guarantee that this is not the case. 

The following proposition will be useful in this section.


\begin{prop} \label{prop:holom} \cite[Thm. 7.4]{Conway}
Let $g:\Omega_0 \rightarrow \mathbb C$, with $\Omega_0 \subset \mathbb C$ open, be holomorphic and $k\in \mathbb R$.
If $G(s_0)=\frac{1}{k}$ for some $s_0\in \Omega_0$ and $m$ be the order of zero  which the function $G(s)-\frac{1}{k}$ has at $s_0$.
Then there exists for every sufficiently small $\varepsilon>0$ a neighbourhood $U_\varepsilon$ of $s_0$ such that the function $G(s)|_{U_\varepsilon}$ attains every value $w$ with $0<|w-\frac{1}{k}|<\varepsilon$ exactly $m$ times.
\end{prop}


The root locus for the class of  infinite-dimensional functions considered here is well-defined.
\begin{thm}\label{thm:rootlocus}
Consider the root locus functions $f_n (k)$ for a minimal system $\Sigma(A,B,C)$.  Then the following statements hold:
\begin{enumerate}
\item  For each $n\in \mathbb N$: $f_n: [0, \infty)\rightarrow \mathbb C$ is well-defined and continuous.
\item  For each $n\in \mathbb N$: $f_n$ is a simple non-intersecting curve.
\item   For each $n\in \mathbb N$: Either there exists an transmission zero $s$ of the system $\Sigma(A,B,C)$ such that $\lim_{k\rightarrow \infty} f_n(k)=s$ or $\lim_{k\rightarrow \infty}| f_n(k)|=\infty$.
\item Let $z \in \rho (A) $ be a transmission zero of the system $\Sigma(A,B,C)$. Then there exists $s_k \rightarrow z$ as $k \rightarrow \infty$ such that $s_k \in \sigma (A-k  BC)$. 
\item For any point $s\in \mathbb C $, only finitely many $f_n$ intersect. Multiplicity of the  spectrum is preserved at such intersection points. 
 Furthermore,  different branches of the root locus do not overlap on any interval.  
\end{enumerate}
\end{thm}

\pf
First prove part 1.
Choose any finite set of eigenvalues of $A$ and  enclose them  by  a simple closed curve $\Gamma$ separating  this part of the spectrum  $\sigma_1 (A)$ from the remainder.  Let $\cal{N}$ be the indices of the eigenvalues of $A$ contained in $\Gamma$. Then \cite[IV.3.5, pg.~213]{Kato} implies there is a $k^M$ such that $f_n (k)$ is a continuous well-defined curve for all $k \in [0 , k^M]$, $n \in {\cal N}$.   Thus, for each $n$,  there is $k^M_n$ so that $f_n$ is a continuous function of $k$ for $k \in [0,k^M_n]$.  
Proposition \ref{prop:holom} together with 
 Lemma \ref{lem:Gk}  implies that the root locus curves are defined on the interval $[0,\infty)$. 
The continuity of the root locus curves follows now from \cite[IV.5, pg. 213]{Kato}.

Part 2 follows directly from Theorem \ref{thm:allk}.  

To show part 3,  assume that $|f_n(k)|$ does not converge to $\infty$ as $k\rightarrow \infty$. 
This implies that there exists
 a sequence $\{y_l\}_{l\in \mathbb N}$  $y_l\ge 0$,  $y_l\rightarrow \infty$, such that $\sup_{l\in \mathbb N}|f_n(y_l)|<\infty$. Due to the Bolzano-Weierstrass Theorem there exists a convergent subsequence, also denoted by $\{ y_l\}, $ such that the limit $f_y:=\lim_{l\rightarrow \infty} f_n(y_l)$ exists. Lemma \ref{lem:Gk} implies $G( f_n(y_{l}))= -\frac{1}{y_{l}}$. As each $\lambda \in \sigma(A)$ is a pole of $G(s)$, see Proposition \ref{prop:mero}, $f_y\not\in \sigma(A)$. Because $G$ is a holomorphic function on $\mathbb C\backslash \sigma(A)$,  $G(f_y)=0$.  Thus every convergent subsequence of $\{f_n(k)\}_{k\ge 0}$ converges to a transmission zero of $G$. Without loss of generality,  assume $f_y=0$. (Otherwise, consider $\Sigma(A-f_y I , B, C)$ which has zero $0$.)

Let $\{ z_l \}$  be another such convergent subsequence with $f_z:=\lim_{l\rightarrow \infty} f_n(z_l)$ and assume $f_z\not=0$.  Define the function $h:[0,\infty)\rightarrow \mathbb R$ by $h(s):=|f_n(s)|$. Due to the continuity of $h$ and the Intermediate Value Theorem for every $l,m\in \mathbb N$,  there exists $x_{l,m}\ge 0$,  such that $h(x_{l,m})= \frac{1}{m} h(z_l)+ (1- \frac{1}{m}) h(y_l)$. Then $\lim_{l\rightarrow \infty} h(x_{l,m}) = \frac{1}{m} |f_z| $. For every $m\in \mathbb N$ choose a subsequence of $\{x_{l,m}\}_l$, which we denote again by $\{x_{l,m}\}_l$ such that $\lim_{l\rightarrow \infty} f_n(x_{l,m}) =: q_m$ exists with $|q_m|=\frac{1}{m} |f_z|  $. Lemma \ref{lem:Gk} implies $G( f_n(x_{l,m}))= -\frac{1}{x_{l,m}}$. As each $\lambda \in \sigma(A)$ is a pole of $G(s)$, see Proposition \ref{prop:mero}, $q_m\not\in \sigma(A)$. Because $G$ is a holomorphic function on $\mathbb C\backslash \sigma(A)$,  $G(q_m)=0$.
Thus the sequence $\{q_m\}$ converges to $0$ and $G(q_m)=0$ for every $m\in \mathbb N$. The fact that  $G$ is holomorphic then implies that $G(s)=0$ on $\mathbb C\backslash \sigma(A)$, which is in contradiction to Proposition \ref{prop:mero}. Thus $f_z=f_y=0 .$
This completes the proof of part 3. 

Next, prove part 4. By the Open Mapping Theorem, point $z \in \rho (A)$ satisfies $G(z ) =0$ if and only if there is $\{ s_n \} \in \mc \cap \rho (A)$ and a real-valued sequence such that $s_n \rightarrow s$, $k_n \rightarrow \infty$ and $G(s_n )= C (s_n I - A)^{-1} B = -\frac{1}{k_n}$.  This is (trivially) equivalent to $ -\frac{1}{k_n} \in \sigma ( C (s_n I - A)^{-1} B) $.  By \cite[pg. 38 (3)]{GGK} this is equivalent to  $ -\frac{1}{k_n} \in \sigma (  (s_n I - A)^{-1} B C )$ and similarly, since $(s_n I -A)^{-1} B $ and $C$ are bounded operators,  $1 \in \sigma (-k_n  (s_n I - A)^{-1} B C) $.  This is also equivalent to $s_n \in \sigma (A-k_n B C)$ 
\cite[Prop. 4.2,pg. 289]{EngelNagel}. 

In order to prove part 5, let $s\in \mathbb C$. If $f_n(k)=s$ for some $n$ and $k$, then $s$ is an eigenvalue of $A-k BC$. As the multiplicity of every spectral point is finite, Proposition \ref{prop:holom} implies that at most finitely many curves $f_n$ intersect at $s$ and the multiplicity of the spectrum is preserved at these intersection points.

To show the last statement, suppose that several branches do overlap on an interval $[a,b]$. 
Since multiplicity is preserved, 
$ G(s) + \frac{1}{k}$ has a zero of at least order $2$ at each $s \in [a,b]$. Differentiating yields $G' (s) =0$ on $[a,b]$. This implies that $G$ is  constant on $\Omega$, which is in contradiction with $G$ being non constant.
\eop

%
%
%

For finite-dimensional systems, where $(A,B,C,D)$ are matrices with real entries, $G(\bar{s}) = \bar{G(s)}$ and so the root locus is symmetric about the real axis. A similar property holds for infinite-dimensional systems, but a definition of a real system is required. 
 \begin{defn}
 If a mapping on $ \cZ $, denoted $x \to \bar{x}$, exists so that for all $x, y \in \cZ,$ and all scalars $\alpha \in \mc$, 
\begin{eqnarray*}
\conj{ (x+y)} &=& \conj{x} + \conj{y}\\
\conj{ \alpha x } &=& \bar{\alpha}\,\, \conj{x} \\
\conj{\conj{x}}&=&x,
\end{eqnarray*}
(where $\bar{\alpha}$ indicates the usual complex conjugate), $\cZ$ is said to have a {\em conjugate operation.}
\end{defn}

\begin{defn}
For any $B \in {\cL} (\cZ_1, \cZ_2)$, where $\cZ_1$, $\cZ_2$ are Hilbert spaces possessing a conjugate operation, $\conj{B} \in {\cL} (\cZ_1, \cZ_2)$ is defined by 
$$\conj{B} z = \conj{ B \conj{z} } .$$
\end{defn}

If $B_1 \in {\cL} (\cZ_1, \cZ_2)$ and $B_2 \in {\cL} (\cZ_2, \cZ_3)$, then it is easy to see that $\conj{B_2B_1}= \conj{B_2}\conj{B_1}$.

\begin{defn}
A  system $\Sigma(A,B,C)$ is called {\em real} if $\conj{A} =A$, $\conj{B}=B$ and  $\conj{C} = C$.
\end{defn}

Although it is easy to write down operators that are not real, control systems arising from partial differential equations and also from delay equations are typically real. 
With these definitions, the following results are straightforward consequences of the above definitions and the definition of a transfer function.

\begin{thm}\label{theoreal}
The transfer function $G$ of a real system $\Sigma(A,B,C)$ satisfies $G(\bar{s})= \bar{G(s)}.$
\end{thm}

\begin{cor}\label{corrlsym}
The root locus  of a real system $\Sigma(A,B,C)$ is symmetric about the real axis. 
\end{cor}

Thus, in summary, the root locus of any control system in the class considered here is well-defined and, if the system is real, symmetric about the real axis. Provided that 
 the system is minimal, then each branch  is a simple, non-intersecting curve. The limit of any branch is a transmission zero or tends to infinity and every transmission zero is the limit of a branch of the root locus.  
 This is similar to the behavior of finite-dimensional systems. 
Unlike the finite-dimensional situation it may happen that the root locus is empty, as the following example shows.

\begin{eg}\label{emptylocus}
The transport equation on the interval $[0,1]$ with a Dirichlet boundary condition is 
\begin{eqnarray*}
\!\!\!\!\!\frac{\partial w}{\partial t}(\zeta, t)\! \!&=&\! \!\frac{\partial w}{\partial \zeta}(\zeta, t)+  bu(t)\quad t\ge 0, \zeta\in[0,1]\\
w(1,t) &=& 0 \quad t\ge 0,\\
w(\zeta,0) \!\!&=&\!\! w_0(\zeta), \quad  \zeta\in[0,1]\\
y(t) \!\!&=& \langle  w(\cdot ,t), c\rangle \qquad t\ge 0.
\end{eqnarray*}
The corresponding operator $A:D(A)\subset L^2(0,1)\rightarrow L^2(0,1) $ is given by $Ax:= x'$ with $D(A)=\{ x\in H^1(0,\infty)\mid x(1)=0\}$ and
$\sigma(A)=\sigma(A-kBC)=\emptyset$.
\end{eg}

In the finite-dimensional case, the number of asymptotes in the root locus (branches converging to infinity) is equal to the relative degree of the system. This is not the case for infinite-dimensional systems. 
The system in the  following example has $\la b, c \ra \neq 0 $ and so is relative degree one but the root locus has  infinitely many branches of the root locus converging to infinity.

\begin{eg}{\bf (Delay Equation)}
Eigenvalues of delay problems are poorly approximated by standard schemes - see for instance, \cite{Silva}.  Furthermore, little is known about zeros or high gain behaviour.   Consider a simple delay equation,
\begin{eqnarray*}
 \dot{x} (t) &=& a x(t) - x(t-1) +u(t), \\ y(t)& =& x(t) . \end{eqnarray*}
 A state-space realization of the form (\ref{sys}) exists on $\cZ= \mc \times L^2 (-1, 0) $ with 
 $$  A \begin{mat} r \\ f(\cdot) \end{mat} = \begin{mat}  a- f(\cdot-1)  \\ f' (\cdot ) \end{mat} , $$
 $$  D(A) = \{ (r, f) \in \cZ;   \;   f  \in H^1 (-1,0),  f(0) =r \} ,$$
 $$B=\begin{mat} 1\\ 0 \end{mat} , \hspace{2em} C= \begin{mat} 1 & 0 \end{mat} . $$
The eigenvalues are given by the roots of 
 \begin{equation} \kappa (s) = s- a +e^{-s} .
 \label{char1}
 \end{equation}
 The invariant zeros are the values of $s$ for which there exist a non-trivial solution $(r, f) \in D(A) $ to the following:
 \begin{eqnarray*}
  s r -  a r + f(\cdot -1 ) + 1& =& 0 ,\\ 
  s f(\cdot) - f' (\cdot) &= & 0 ,\quad
  r = 0.
  \end{eqnarray*}
  The only solution to this system of equations is the trivial solution and so there are no invariant zeros. 
 Since
  $${\rm rank } \begin{mat} \kappa (s) ; 1 \end{mat} = 1,$$ 
the system is approximately controllable \cite[Thm. 4.2.10]{CZbook} and since
 $${\rm rank } \begin{mat} \kappa (s) \\ 1 \end{mat} = 1,$$ 
 the system is approximately observable \cite[Thm. 4.2.6]{CZbook}.
 Since the systems is approximately controllable and observable, these same conclusions can be found by examining the transfer function
  $$  G(s) =  \frac{1}{\kappa (s) }. $$

 The eigenvalues of $A$ form a sequence with $\re \lambda \rightarrow -\infty$ as $|\lambda | \rightarrow  \infty$ and in fact
 $$ | \lambda |  \leq |a | + e^{\re \lambda} .$$
 \cite[Prop. 1.8, Prop. 10]{MichielsNiculescu}.
 The eigenvalues of $A-k B B^*$ are the roots of
\begin{equation}  s- a +k +e^{-s} \label{char} \end{equation}
and so they have a similar pattern, for each $k$, as the eigenvalues of $A$.

\begin{thm} \cite[Thm. 6.1]{Silva}
Consider the  equation
$$ \delta (s) = k_p k_c e^{-s} + 1+ T s$$
where $T>0$, $k_p >0$.  All roots of this equation will have negative real parts if
$$ -\frac{1}{k_p} < k_c< \frac{T}{k_p} \sqrt{z_1^ 2 + \frac{1}{T^2} }$$ 
where $z_1 \in (\frac{\pi}{2}, \pi) $ solves
$$ \tan (z) = - Tz . $$
\end{thm}
\smallskip
Rewriting the characteristic equation (\ref{char}) in the above form reveals that all the eigenvalues of $A+k B B^*$ are stable for every $k> a$.
 Since there are no zeros,  all branches of the root locus move from the eigenvalues of $A$  to $-\infty$  in the left-half-plane. 
\end{eg}

 Now consider the root locus for systems with non-zero feedthrough $\Sigma(A,B,C,D)$, that is,
\begin{equation}
\begin{array}{lll}
\dot{x}(t) &=& Ax(t)+Bu(t)\\
y(t) &=& Cx(t) + Du(t),
\end{array}
\label{sysd}
\end{equation}
where $A,B,C$ are as above but $D\not=0$.

\begin{lem}
The transfer function of the system $\Sigma(A-B\Di C, B \Di , - \Di C , \Di )$ is  $G(\cdot)^{-1}$.
\end{lem}

\pf This is a straightforward calculation since $B$ and $C$ are bounded.  \eop

The following result is straightforward. 
\begin{lem}
Suppose $\Sigma(A,B,C,D)$ is a system with $D\not=0$.
Then for $k>-\Di$, define $A_{k}=A-B k (1+Dk )^{-1} C$ with $D(A_{k}) =D(A)$. The system with feedback $u(t)=-k y(t)+v(t)$ has realization $\Sigma (A_{k} , B(I+k D)^{-1} , (1+Dk)^{-1} C , (1+Dk)^{-1} D)$.
\end{lem}

In a similar manner as for $D=0$,   define the root locus for general systems $\Sigma(A,B,C,D)$. The root locus is  the eigenvalues of  $A-B k (1+Dk )^{-1} C$ for real $k$. Since the set $ \sigma(A-B k (1+Dk )^{-1} C)$ is finite or countable consisting of eigenvalues of $A-B k (1+Dk )^{-1} C$   there is a family of curves $f_n (k)$ associated to the eigenvalues of $A$ with $f_n (0)= \lambda_n$. The values of $f_n (k)$ are the eigenvalues of  $A-B k (1+Dk )^{-1} C$.   The {\em root locus} is the set of curves $f_n (k)$.

\begin{thm}
For a minimal system $\Sigma(A,B,C,D)$ where $D \neq 0$, the root locus is well-defined for all $k>0$ and non-intersecting. Each branch moves from a pole to a zero.
\end{thm}

\pf The root locus corresponding to $\sigma (A-B k (1+Dk )^{-1} C) $ as $k$ increases from $0$ to infinity  is the same as that of $A-\tilde{k} BC$ where $\tilde{k}$ increases from $0$ to $\Di$.  Thus the root locus is well-defined by Theorem \ref{thm:rootlocus}. But $\sigma (A-\Di BC)$ are the zeros of $\Sigma(A,B,C,D)$ so each branch moves from a pole to a zero.  \eop

\section{Spectrum Determined Growth Assumption}

The relevance of the root locus to control system design and analysis relies on a relationship between the spectrum of $A-kBC$ and the growth of the semigroup it generates, and hence the dynamics of the system. 
For finite-dimensional systems, the spectrum of  $A-Bk(1-Dk)^{-1}C$  determines the dynamics of the controlled system. However, for infinite-dimensional systems, this is not always the case. 
Systems for which the spectrum of the generator does determine the growth (or decay) of the associated semigroup as said to satisfy the following assumption. 
\begin{defn}
Let $A$ be the generator of a $C_0$-semigroup $(S(t))_{t\ge 0}$ on a Hilbert space $\cZ$. The generator is said satisfy the {\em Spectrum Determined Growth Assumption} (SDGA) if $ \sup_{\lambda \in \sigma (A) }  {\rm Re} \lambda  =\inf\{ \omega \in \mathbb R\mid \,\exists \, M>0 : \| S (t) \| \leq M e^{\omega t},\, t\ge 0\}$.
\end{defn}
Unless the  SDGA  holds  for  $A-Bk(1-Dk)^{-1}C$ for all $k \geq 0$, spectral analysis of  the generator $A-Bk(1-Dk)^{-1}C$ says nothing about growth or stability of the semigroup.

The SDGA holds for analytic semigroups,  differentiable semigroups, compact semigroups and  also Riesz spectral systems and delay equations \cite[Cor. 3.12, p. 281]{EngelNagel}, \cite[Thm 2.3.5 and Theorem 5.1.7]{CZbook}. Further the SDGA assumptions hold for a generator $A$ of a $C_0$-semigroup,
satisfying the following conditions (see  \cite[Theorem 3.5]{XuFeng}):
\begin{enumerate}
\item[(A1)] $\sigma(A)$ consists of isolated eigenvalues with finite algebraic multiplicity;
\item[(A2)] all eigenvalues of $A$ with sufficiently large module are semi-simple, that is, its geometrical multiplicity is equal
to the algebraic multiplicity;
\item[(A2)]  the sequence of the generalized eigenvectors of $A$ forms a Riesz basis in ${\cal Z}$;
\item[(A4)] for any $r > 0$, there exists an integer $M$ such that, for any given $s\in \mathbb C$, we have
\[ {\rm dim}\, \left( \sum_{ \lambda \in \sigma (A)\cap B(s,r)} P_\lambda {\cal Z}\right) \le M,\]
where $P_\lambda {\cal Z}$ is the spectral projector corresponding to the eigenvalue $\lambda$.
\end{enumerate}
The SDGA assumption is preserved under bounded perturbations such as $Bk(1-Dk)^{-1}C$ for analytic semigroups, compact semigroups and $C_0$-semigroups satisfying (A1)-(A4), see \cite[Thm. 4.3]{pa}, \cite{NagelPiazzera} and \cite[Theorem 3.5]{XuFeng}.

The SDGA is not in general preserved under bounded perturbations of a generator of a differentiable semigroup; see the counter-example on a Hilbert space in \cite{Renardy}. 
However,  a positive result can be obtained for  the perturbations that arise in  control of delay equations.
Consider a general class  of delay equations, see \cite[Section 2.4]{CZbook},
\begin{eqnarray}
\dot{x}(t) &=& A_0 x(t) + \sum_{i=1}^p A_i x(t-h_i) + b_0 u(t), \qquad t\ge 0,\nonumber\\
x(0) &=& r,\label{eqn:delay}\\
x(\theta) &=& f(\theta), \qquad  -h_p\le \theta<0,\nonumber
\end{eqnarray}
where $0< h_1<\cdots <h_p$ represents the point delays, $x(t)\in \mathbb C^n$, $A_i\in \mathbb C^{n\times n}$, $i=1, \cdots, p$, $r\in \mathbb C^r$, $b_0\in \mathbb C^n$ and $f\in L^2 (-h_p,0;\mathbb C^n)$.
Equation \eqref{eqn:delay} can be reformulated as an abstract differential equation of the form
\[ \dot{z}(t) =A z(t) + bu(t), \qquad z(0)=z_0,\]
where $A$ generates a $C_0$-semigroup on $\cZ:=\mathbb C^n \times L^2 (-h_p,0;\mathbb C^n)$ and $b\in \cZ$. In \cite[Theorem 5.1.7]{CZbook} it is shown that $A$ satisfies the SDGA.

\begin{thm}\cite[Theorem 5.1.7]{CZbook}
Consider the delay equation \eqref{eqn:delay} with measurement
\[ y(t) = c_0 x(t),\]
where $c_0\in \mathbb C^{1\times n}$ and the feedback control
\[ u(t) = -k y(t)+v(t)\]
where $v$ is an exogenous input. 
Then the  closed loop system is again a delay equation of the form \eqref{eqn:delay}   with $A_0$ replaced by $A_0 -k b_0c_0$ and the corresponding generator $A-k BC$ satisfies the SDGA.
\end{thm}
Thus, the root locus is a useful tool for analyzing stability of the controlled  infinite-dimensional systems that commonly arise in applications. 

\section{Collocated Self-Adjoint Systems}


Many diffusion problems, such as heat flow, lead to a system where the generator $A$ is self-adjoint. 

\begin{defn}
A system $\Sigma(A,B,C)$ is called {\em collocated self-adjoint}, if $A$ is self-adjoint and  $C=B^*$.
\end{defn}
 The self-adjoint operator $A$ is {\em negative semi-definite}, if $ \langle Az, z\rangle\le 0$ for every $z\in D(A)$.

If the underlying state space is finite-dimensional, then it is well-known that the poles and zeros are real, interlace on the negative real axis and furthermore, the system is relative degree one so that there is one asymptote. This asymptote moves along the negative real axis to $-\infty$.  A partial generalization for infinite-dimensional systems  was obtained in \cite{ZwartHof}. In that paper the authors show that the invariant and transmission zeros are real and that the poles and zeros  of Riesz spectral systems that satisfy an additional technical condition  interlace on the real axis. 
The following theorems, which use results provided earlier in this paper,  provide  a significant generalization of this earlier work.

\begin{thm}
 \label{col}
Suppose that the system $\Sigma(A,B,B^*)$ is minimal and collocated self-adjoint. 
Then each $f_n (s)$ is real-valued, all the zeros are real and interlace with the eigenvalues.  If the state space is infinite-dimensional,  the root locus has no asymptote and each branch converges to 
the left to a zero as $k \rightarrow \infty$.
 \end{thm}
\pf
It is well-known that all the eigenvalues of $A$ are real and since $A-k B B^*$ is also self-adjoint
all branches of the root locus lie entirely on the
real axis.   Since $A$ is assumed to be a generator, the  real part of the eigenvalues is bounded (by the growth of the semigroup) and the eigenvalues form a sequence moving to $- \infty$ on the real axis. 

As shown for general systems in section III,  no branch intersects with itself and so each branch starts at an eigenvalue and moves monotonically either to the left or right.
Minimality of the system implies that the root locus does not intersect with $\sigma (A)$ for any value of $k$. 

Each branch of any root locus  converges either to a zero or to  infinity. It will now be shown that each branch moves to the left from an eigenvalue and,  for infinite-dimensional systems,  converges to a zero.
Let $x_n(k)$ be a normalized eigenvector corresponding to the eigenvalue $f_n(k)$ of $A-kBB^*$, that is, $$(A-kBB^*)x_n(k)=f_n(k)x_n(k).$$
This is a differentiable function of $k$ \cite[Lemma 4.7]{Langer} and so
\[ -BB^*x_n(k)+ (A-kBB^*)x_n'(k)=f_n'(k)x_n(k)+f_n(k)x_n'(k).\]
Taking the inner product with $x_n(k)$,
\[ -|B^* x_n(k)|^2 = f_n'(k).\]
Therefore, each branch moves to the left.   

If the space $\cZ$ is infinite-dimensional, the spectrum of $A$ consists out of an unbounded sequence $\{\lambda_n\}$ of negative real numbers.  Since the root locus cannot intersect with an eigenvalue, each branch must converge. But any bounded branch must converge to a zero. Thus each branch moves to the left to a zero as $k \rightarrow \infty$, and the root locus has no asymptote. This also implies that the zeros interlace with the eigenvalues. 
\eop

If instead  $k \rightarrow -\infty$, the above argument yields that each branch
of the root locus converges to the right and that there could be one branch of the root locus that converges to $ \infty$.

If $A$ is defined on a finite-dimensional space, then there are a finite number of eigenvalues and there is always one branch of the root locus that converges to $- \infty$ as $k \rightarrow \infty.$ 

The following result shows that  zeros of $A$  become arbitrarily close to  the eigenvalues for large $|s|.$
 \begin{thm}
 \label{zeros_poles}
 Consider the minimal, collocated self-adjoint system $\Sigma(A,B, B^*)$ with  $A$ negative semi-definite and  $b \in D((-A)^{\half})$. Then, indicating the invariant zeros by $\{ \mu_n\}, $ for any $\eps>0$ there is $N\in \mathbb N$ 
 so that for every $n\ge N$, 
 $|\mu_n - \lambda| < \eps$ for some eigenvalue $\lambda $ of $A.$
 \end{thm}

 \pf  The invariant zeros of $\Sigma(A,B,B^*)$ and  $\Sigma(A, (-A)^{\half} B  , B^* (-A)^{-\half} )$ are identical, see Lemma \ref{lem:alpha}.
Since $A$ is self-adjoint, the operator $A_\infty$, given by \eqref{Ainfinity},
\[ A_\infty z = Az-\frac{1}{B^*B} (-A)^{\half}BB^*(-A)^{\half},\quad z\in D(A),\]
is self-adjoint and therefore the eigenfunctions $\{z_n \}  \subset \kerc $ of $A_\infty$ are orthogonal and they can of course be chosen normal.
Applying Theorem \ref{zeros_close} to the system $\Sigma(A, (-A)^{\half} B  , B^* (-A)^{-\half} ) $ leads to the conclusion that $\{ \mu_n \}$ are
in the $\eps$-pseudo\-spectrum of $A$. Since  $A$ is self-adjoint this means that   
for any $\eps>0$ there is $N$ so that  for all $n\ge N$, there is $\lambda  \in \sigma (A)$ such that  $| \lambda - \mu_n | < \eps. $
 \eop

 Since $A_k = A-k B B^* = A $ on ker$\, B^*$, it follows by a similar calculation  that  the invariant zeros are in the $\varepsilon$-pseudospectrum of the root locus, uniformly in $k$.
 Theorem \ref{zeros_poles} generalizes  \cite[Thm. 4.4]{ZwartHof}  where additionally assumptions on the eigenvalues are required.  However,   convergence rates are obtained in this earlier work. \\

The above results are illustrated by the following simple example. 

\begin{eg}{\bf (Heat flow in a rod)}
Consider the problem of
controlling the temperature profile in a rod of length $1$ with constant thermal conductivity $\kappa$, 
mass density $\rho$ and specific heat $C_p$. The rod is insulated at the ends $x=0$, $x=1$. To simplify, use dimensionless variables so that $\frac{\kappa}{C_p \rho}=1$.
With control applied through some weight $b(x)$, 
  and the temperature is governed by the following problem
$$\frac{\partial z(x,t)}{\partial t} = 
   \frac{\partial^2 z(x,t)}{\partial x^2} +b(x) u(t),   \qquad x\in
      (0,1), \qquad t \geq 0. 
      $$
      $$
         \frac{\partial z}{\partial x}(0,t) = 0,      
      \qquad \frac{\partial z}{\partial x}(1,t) =  0,
      $$
      where $b (x) \in L^2 (0, 1)$. 
  The temperature sensor is modelled by
  \begin{equation} 
  y(t) = \int_0^1 b(x) z(x) dx .
  \label{heatc}
  \end{equation}
  It is well-known that this can be written as an abstract control system (\ref{sys}) on the Hilbert space $L^2 (0,1)$ with 
  $$ Az =  \frac{\partial^2 z}{\partial x^2} , \hspace{2em}  D(A) = \{ z \in H^2 (0,1) , \; z' (0 ) = z' (1) = 0 \} $$
  and $Bu =b (x) u$, and $C=B^*$ is defined by (\ref{heatc}).  
  This system is approximately controllable (and observable) if 
  \begin{equation}
   \int_0^1 b (x)  \cos (n \pi x ) dx \neq 0 
   \label{bcont}
   \end{equation}
    for all integers $n$ \cite[Thm. 4.2.1]{CZbook}. Since the eigenfunctions form an orthonormal basis for $L^2 (0,1)$, the control system is minimal if (\ref{bcont}) is satisfied.
    The operator $A$ is a self-adjoint, negative semi-definite  operator. Thus, the eigenvalues $\lambda$ of $A$ are all real and non-positive. In this case,  $\lambda= - n^2 \pi^2$, $n \geq 0$. 
    
    The invariant zeros depend on $b(x)$, but Theorem \ref{col} implies that they are real and negative (since 0 is an eigenvalue of $A$, it cannot be a zero) and if $b(x)$ satisfies (\ref{bcont}), the zeros interlace with the eigenvalues and the  eigenvalues of $A-kB B^*$ converge to the zeros.  
    
    The system with simple proportional control $u=-ky+v$ is always stable, but its decay rate is limited by the largest zero, which lies in the interval $[-\pi^2 , 0]$.
    The largest invariant zero is found by solving
 $$ s w_0 - w_0^{\prime \prime } = b $$
 for $w_0 (x,s)$ and then finding the root of 
 $$ p(s)= \int_0^1 w_0 (x,s) b(x) \, dx   $$
 that lies in $(-\pi^2 , 0 ). $
 If  $b(x)=x^2,$ this largest zero is $ -5.65.$  The next two zeros are $-38.6$ and $-88.4, $ quite close to the eigenvalues $-39.5$ and $-88.8$ respectively.
 The largest zero, -5.65, determines the limit of the settling time that can be achieved with constant gain feedback.
The qualitative nature of the root locus for a collocated self-adjoint system  is illustrated in Figure \ref{fig3}.
\begin{figure}
\centering
\includegraphics[width=3.5in,height=2.in]{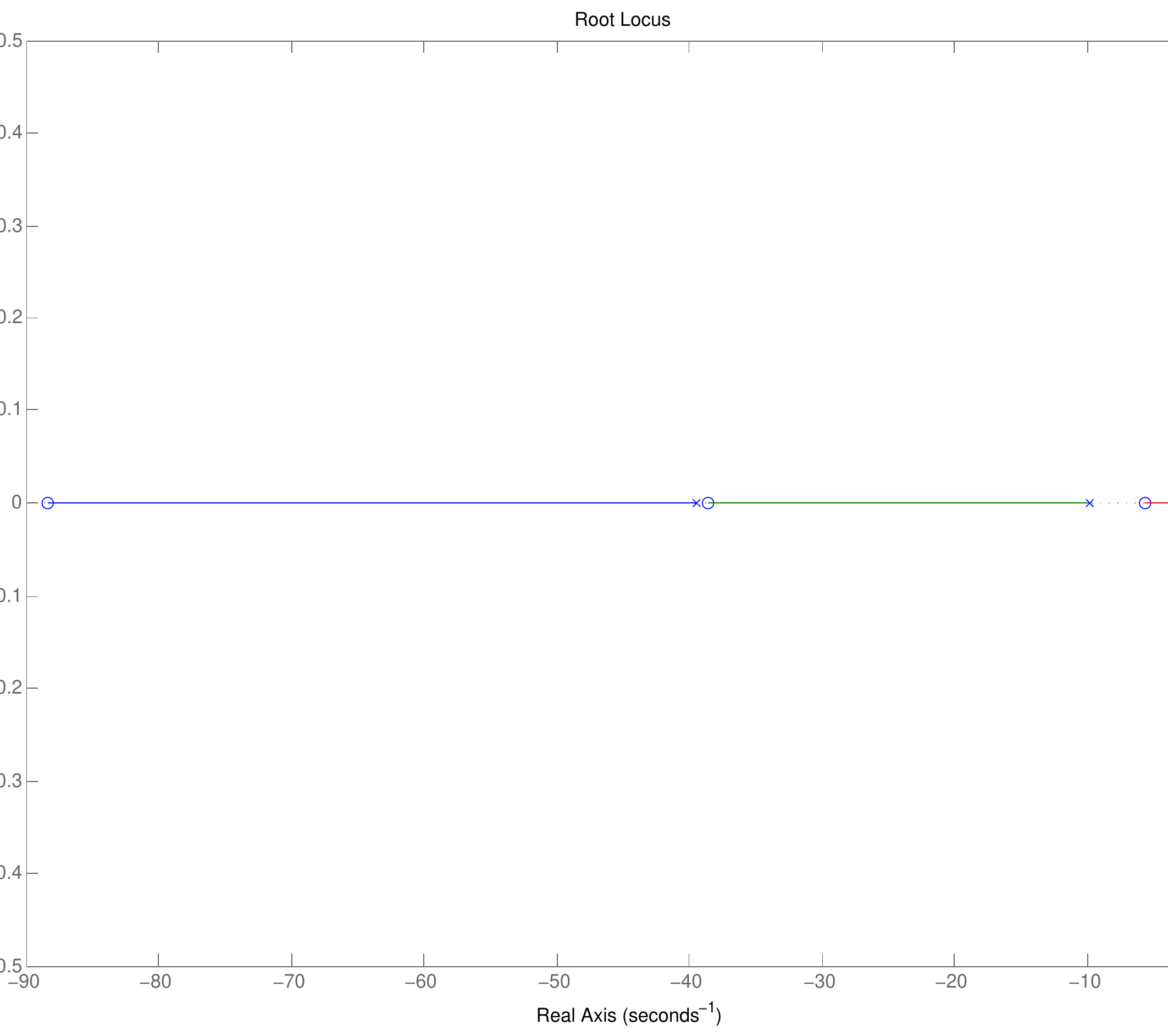}
\caption{The root locus for a collocated self-adjoint system lies entirely on the real axis. For an infinite-dimensional system, each  branch starts at a pole and moves to the left towards a zero. This figure was calculated using the first 3 eigenvalues and first 3 zeros of the heat flow example. \label{fig3} }
\end{figure}

\end{eg}

 \section{Collocated Skew-Adjoint Systems}

Undamped waves and structural vibrations lead to a control system where the generator satisfies $A^*=-A;$ that is the generator is skew-adjoint.

\begin{defn}
A system $\Sigma(A,B,C)$ is called {\em collocated skew-adjoint}, if $A$ is skew-adjoint and  $C=B^*$.
\end{defn}

  In \cite[Thm. 3.2]{ZwartHof} it is shown that the transmission zeros lie in the left half-plane. Here it is  shown that the zeros interlace with the eigenvalues on the imaginary axis and that the root locus lies in the left half-plane. Furthermore, under an additional assumption, the zeros become asymptotically close to the eigenvalues.

\begin{thm}\label{theoskewinterlace}
If $\Sigma(A,B,B^{*})$ is a minimal collocated skew-adjoint system 
 then the  zeros and eigenvalues  lie on the imaginary axis and the zeros interlace the eigenvalues. 

Moreover, each branch of the root locus lies in the closed left half plane. More precisely, if $\lambda \in \sigma (A-k B B^{*})$, $k>0$ then $\re \lambda <0 $.  
\end{thm}

\pf  
Since $BB^*$ is a compact operator, the spectrum of  $A-k B B^{*}$ consists of isolated eigenvalue only  for every $k>0$.
Thus, if $ \lambda   \in \sigma (A-k B B^{*})$,  there is  some  $x \in \cZ$, $\|x\| =1$, so that
$$ \lambda x = A x - k B B^* x .$$
This implies
\begin{equation*} \lambda= \la Ax, x \ra - k | B^* x |^2 \quad\mbox{and} \quad
\bar{\lambda} = \la x, Ax\ra -  k | B^* x|^2 . \label{eig:skew} \end{equation*}
Since $A$ is skew-adjoint, 
$$ 2\re \lambda = -k | B^* x|^2 .$$
This is always non-positive, and negative if $k>0$ since the system is minimal.  

Define $\tilde A:=iA$ with $D(\tilde A)=D(A)$. 
Then $\tilde A$ is a self-adjoint operator, but in general not the generator of a $C_0$-semigroup. However, the proof of Theorem \ref{col} implies that the transmission and invariant zeros of $\Sigma(iA,B,B^*)$ are real and the zeros interlace with the poles if  the system $\Sigma(iA,B,B^*)$ is minimal. 
 
 An easy calculation shows that $\lambda$ is a transmission zero, invariant zero or eigenvalue of $\Sigma(iA,B,B^*)$ if and only if $-i\lambda$ is a transmission zero, invariant zero or eigenvalue of $\Sigma(A,B,B^*)$, respectively. Thus, all the zeros are imaginary. Also, the minimality of $\Sigma(A,B,B^*)$ implies that the transmission and invariant zeros interlace.
 \eop

\begin{prop}
Consider a real, minimal, collocated skew-adjoint  system $\Sigma(A,B,B^{*})$. Then zero is either an invariant zero
or $0 \in \sigma (A)$.
\end{prop}

\pf This follows from the fact that the eigenvalues  and zeros  interlace on the imaginary axis, see Theorem \ref{theoskewinterlace}, and that the root locus is symmetric about the real axis, see Corollary \ref{corrlsym}.
\eop

 \begin{prop}
Consider a real, minimal,   collocated skew-adjoint  system $\Sigma(A,B,B^{*})$.
 Then the entire negative real axis is part of the root locus.
 \end{prop}
 
 \pf
 Since $A$ is skew-adjoint, $A$ generates a unitary $C_0$-group and  $G(s) = B^* (sI -A)^{-1} B $ is an analytic function for $\re s <0$. Since the system is real, $G(\bar{s}) = \bar{G(s)}$  (Theorem \ref{theoreal}) and so $G(s) \in \mathbb R$ if $s<0. $
  For any $s<0$ define $x_s:= (sI-A)^{-1}B$. Then since not only is the system real, but $A^*=-A,$
\begin{eqnarray*}
G(s)&=& B^* (sI-A)^{-1}B \\
&=& x_s^*(sI+A)x_s \\
 &=& s \|x_s\|^2 + x_s^*Ax_s  \\
 &=&  s \|x_s\|^2\\
 &< & 0.
\end{eqnarray*}
 Thus, for any $s<0$ there is $k>0$ so $G(s)= -\frac{1}{k}$ and
 Lemma \ref{lem:Gk}  implies that the entire real axis is in the root locus.
  \eop

  \begin{thm}
Consider a  minimal, collocated skew-adjoint  system $\Sigma(A,B,B^{*})$.
Then the  semigroup generated by $A-kBB^*$, $k>0$, is strongly stable. 
  \end{thm}
  
 \pf As $A$ is skew-adjoint, $A$ generates a unitary $C_0$-group. Moreover, $-kBB^*$ is a bounded dissipative operator on ${\cal Z}$ for $k>0$. This implies that  $A-kBB^*$, $k>0$, generates a contraction semigroup, see \cite[Chapter III, Theorem 2.3]{EngelNagel}. Theorem \ref{theoskewinterlace} implies that the spectrum of $A-kBB^*$, $k>0$, is contained in the open left half plane. Now the statement of the theorem follows from the Arendt-Batty-Lyubich-V${\rm \tilde{u}}$-Theorem, see \cite[Theorem V.2.21]{EngelNagel}.
 \eop

 \begin{thm} \label{theo:colconv}
 Consider the mimimal collocated skew-adjoint system $\Sigma(A,B, B^*)$ with \,
$b \in D((-A)^{\half})$ where $Bu=bu$. Indicating the invariant zeros by $\{ \mu_n\}, $ for any $\eps>0$ there is $N\in \mathbb N$ so that for every $n\ge N$, 
 $|\mu_n - \lambda| < \eps$ for some eigenvalue $\lambda $ of $A.$
 \end{thm}
 
 \pf The invariant zeros of $\Sigma(A,B,B^*)$ and  $\Sigma(A, (-A)^{\half} B  , B^* (-A)^{-\half} )$ are identical, see Lemma \ref{lem:alpha}.
A straightforward calculation shows that the  eigenfunction of  $A_{\infty} $  (defined in Theorem \ref{bcne0}) corresponding to the eigenvalue $\mu_n$ is $ (\mu_n I-A)^{-1}b$. Using the resolvent identity and the fact that $\mu_n\in i\mathbb R$,
yields that  for $n\not=m$
\begin{align*}
\langle ( & \mu_n I-A)^{-1}b,   (\mu_m I-A)^{-1}b\rangle  \\ & = -\langle  (\mu_m I-A)^{-1}(\mu_n I-A)^{-1}b, b\rangle\\
&= \frac{1}{\mu_n-\mu_m} ( \langle  (\mu_n I-A)^{-1}b, b\rangle -\langle  (\mu_m I-A)^{-1}b, b\rangle)\\ 
&=0,
\end{align*}
because $\mu_n$ and $\mu_m$ are invariant zeros of $A$. 
Let $$z_n:= \| (\mu_n I-A)^{-1}b\|^{-1}  (\mu_n I-A)^{-1}b. $$ Then the set
 $\{z_n \}$   is orthonormal. 
Applying Theorem \ref{zeros_close} to the system $\Sigma(A, (-A)^{\half} B  , B^* (-A)^{-\half} ) $ leads to the conclusion that $\{ \mu_n \}$ are
in the $\eps$-pseudo\-spectrum of $A$. Since  $A$ is normal this means that   
for any $\eps>0$ there is $N$ so that  for all $n\ge N$, there is $\lambda  \in \sigma (A)$ such that  $| \lambda - \mu_n | < \eps. $
 \eop
 
 Since $A_k = A-k B B^* = A $ on ker$\, B^*$, it also follows that  the invariant zeros are in the pseudospectrum of the root locus, uniformly in $k$.
 
 \begin{cor}
  Consider the minimal collocated skew-adjoint system $\Sigma(A,B, B^*)$ with \,
$b \in D((-A)^{\half})$ where $Bu=bu$. Indicating the invariant zeros by $\{ \mu_n\}, $ for any $\eps>0$, there is $N$ so that for all $k$ and 
 for all $n>N$,  $\mu_n \in \sigma_\eps(A_k)$.
 \end{cor}
 
 For a special class of skew-adjoint systems that occurs often in applications, second-order systems, a rate of convergence of the zeros to the eigenvalues can be obtained. 
First define on a Hilbert space $H$ the stiffness operator $A_o : D (A_o) \subset H
\rightarrow H$ to be a self-adjoint, positive-definite linear operator
 such that zero is in the resolvent set of $A_o$. Here
$D(A_o)$ denotes the domain of $A_o$.
 Since $A_o$
is self-adjoint and positive-definite,  $A_o^\half$ is well-defined. The Hilbert space $V$ is defined
to be $V=
D(A_o^\half) $ with the norm induced by 
\[ \langle x,z\rangle_{V}=
\langle A_o^\half x ,A_o^\half z\rangle_H,\qquad x,z \in V.\]
 Define then, for  $F \in \cl (\mc, H) $, a class of second-order systems
 \begin{equation}
 A = \begin{mat} 0 &I \\ -A_o & 0 \end{mat} , \quad B = \begin{mat} 0 \\ F \end{mat} .
 \label{2nd}
 \end{equation}
 Let $f \in H$ indicate the element of $H$ that defines $F.$
This class describes, for instance, undamped wave and structural vibrations.
It is well-known that $A$ with
domain 
$D(A)= D(A_o)\times V  $
is a skew-adjoint operator on $\cZ= V \times H$ and generates a unitary semigroup on $\cZ .$

\begin{thm}
Consider the system $\Sigma(A,B,B^*)$ where $A,$ $B$ are defined in (\ref{2nd}) and assume that $f \in D(A_o^{\frac{3}{8}} )$ and $Fu=fu$.  Further, we assume that the system is minimal.
Then, indicating the invariant zeros by $\{ \mu_n\}, $ there exists a constant $M>0$ such that  for every $n\in \mathbb N$ 
$$|\mu_n - \lambda| < \frac{M}{\sqrt{|\mu_n|}}$$ for some eigenvalue $\lambda $ of $A.$
%
\end{thm}

\pf
Consider the system  $\Sigma(A,\tilde{B}, \tilde{C})$ where 
$$\tilde{B}=\begin{mat} 0 \\ A_o^{\frac{3}{8}} f \end{mat} , \quad \tilde{C}\begin{mat} w\\ v\end{mat}  = \la v , A_o^{-\frac{3}{8}} f \ra .$$
As shown in Lemma \ref{lem:alpha}, this system has the same zeros as $\Sigma(A,B,B^*).$  Indicate the  eigenfunctions of $A_\infty$ corresponding to $\{ \mu_n\}$ by $\{ z_n  \}$. It is easy to see that the normalized eigenfunctions are of the form  $z_n =[ w_n , \mu_n w_n] ,$ where $w_n\in D(A_o)$. Note that
\[   A_\infty \begin{mat} w\\ v\end{mat} = A \begin{mat} w\\ v\end{mat} - \begin{mat} 0 \\   \frac{\langle w,  A_o^{-\frac{3}{8}}  f\rangle_V}{\|f\|^2}  A_o^{\frac{3}{8}}  f\end{mat}   \]
Therefore,
\begin{eqnarray}
\| \mu_n z_n - Az_n \| &=&\| A_\infty z_n - Az_n \| \\ &=& \frac{1}{\|f\|^2} \|\la A_o w_n , A_o^{-\frac{3}{8}}f \ra  A_o^{\frac{3}{8}} f \| \nonumber\\
&\leq &  \frac{\|A_o^{\frac{3}{8}} f \|}{\|f\|^2} |\la A_o^{\frac{1}{4}} w_n , A_o^{\frac{3}{8}} f \ra  |  \nonumber\\
&\leq &  M\sqrt{2} \|A_o^{\frac{1}{4}} w_n  \|, \label{bd}
\end{eqnarray}
where $M= \frac{\|A_o^{\frac{3}{8}} f \|^2}{ \sqrt{2}\|f\|^2}$.
Moreover, 
\begin{eqnarray*}
1 =\| z_n \|^2 &=& \|A_o^{\half} w_n \|^2 + |\mu_n|^2 \| w_n\|^2  \\
&\geq & 2 |\la A_o^{\half} w_n, \mu_n w_n \ra | \\
&= & 2 |\mu_n | \| A_o^{\frac{1}{4} } w_n \|^2 
\end{eqnarray*}
and 
\begin{equation}
\|A_o^{\frac{1}{4} } w_n\|^2 \leq \frac{1}{ 2 |\mu_n|} .
\end{equation}
 Substituting this bound into (\ref{bd}) yields 
$$\| \mu_n z_n - A z_n \| \leq \frac{M}{\sqrt{ |\mu_n|}}.$$
Since $A$ is skew-adjoint, and hence normal,  the result follows. 
\eop

The above results are illustrated by the following examples.

\begin{eg}
 {\bf (Wave equation on an interval)}
The wave equation  models vibrating strings and many other situations 
such as acoustic plane waves, lateral vibrations in beams, and electrical transmission lines. 
Suppose that the ends are fixed with control and observation  both distributed along the string. For simplicity normalize the units 
 to obtain the equations
\begin{eqnarray}
\label{notdamped}
\frac{\partial^2 w}{\partial t^2}&=&\frac{\partial^2 w}{\partial x^2}+f(x)u(t) , \\  
w (0,t)&=&0, \qquad w(1,t)=0,  \label{bc} \\
y(t)&=& \int^1_0 \frac{\partial w}{\partial t} (x,t) f(x){ d}x , \label{obs}
\end{eqnarray}
where $f\in L^2(0,1)$ describes both the actuator and sensing devices.
The zeros can be found from the zeros of the transfer function, or equivalently from calculating the invariant zeros.  These are the values of $s$ for which
\begin{equation}
\begin{array}{lll}
s^2 w_0(x) -w_0^{\prime\prime} (x) &=& f(x) , \quad w_0(0)=0, \, w_0(1) =0 
\\
s \int_0^1 w_0 (x) f(x) dx &=& 0 
\end{array}
\label{eq:invzero}
\end{equation}
has a non-trivial solution for $w_0 .$
Calculation of the root locus is similar to solving (\ref{eq:invzero}).  Since $u= - k y=-k \langle v, f\rangle , $
\begin{eqnarray*}
s w &= &v\\
sv &=&w^{\prime\prime} - f k \langle v, f \rangle
\end{eqnarray*}
or
$$ s^2 w =w^{\prime\prime} - f s k \langle w, f \rangle. $$
This is the same problem as (\ref{eq:invzero}), except, defining  $\alpha=\langle  w, f \rangle ,$
$f$ is changed to $(-sk \alpha) f . $
Since the problem is linear
$w_k (x)= -sk \alpha w_0 (x) .$
But it is also required that 
$$\alpha=\langle w_k ,f \rangle = -s k \alpha \langle w_0 , f \rangle. $$
If $\alpha=0$, $s$ is an invariant zero. If $\alpha \neq 0$ divide through by $\alpha$ and rearrange to obtain
$$\frac{1}{k}=-s \langle w_0 , f\rangle.$$
The root locus is found by finding solutions to this equation for each $k.$  Calculating the transfer function $G$ and solving $G(s)=-\frac{1}{k}$ yields an identical calculation. 

Setting $f(x)=x$ yields, since the eigenfunctions are $\begin{bmatrix} \phi_n \\ \lambda_n \phi_n \end{bmatrix} $ where $\lambda_n =-\imath n \pi $ and  $\phi_n (x) = \sin (n \pi x ) $ ,
$$ \langle f , \phi_n \rangle = \imath cos(n \pi ) \neq 0 .$$
The system is approximately  controllable and observable and
the transfer function  is
\begin{eqnarray*}
G(s) &=&\frac{ s^2 (1-e^{-2s } ) -3 s (1+e^{-2s } ) +3 (1-e^{-2s } )}{3s^3 (1-e^{-2s } )}\\
&=&\frac{ (s^2+3s+3) (1-e^{-2s } )  -6s }{3s^3 (1-e^{-2s } )}
 . 
\end{eqnarray*}
The first 4 eigenvalues are 
$$ 3.14 , \; 6.28 , \; 9.42 , \; 12. 6 $$
while the first 4 zeros are
$$0, \; 5.76 , \; 9.09, \; 12.3$$
showing that the zeros rapidly become quite close to the eigenvalues. This suggests that in numerical calculation of zeros, the eigenvalues can be used as initial estimates in an iterative algorithm for finding generalized eigenvalues.

The root locus is found by solving 
$$\frac{1}{k} = - G(s)$$
for $s, $ or 
$$\frac{s}{k}=\frac{-1}{3}+\frac{1}{s} (-1 +{2}{(1-e^{-2s } )}   -\frac{1}{s^2} .
$$
Recall that the entire real axis is included in the root locus and also that each zero is a limit of a branch of the root locus. 
If $|s| \to \infty$ and $1-e^{-2s}$ is bounded below then $s$ becomes real and in fact $s \approx -\frac{k}{3} .$  Alternatively, $1-e^{-2s} \to 0$ which means that the root locus is approaching a pole which will be near zero. 
Thus, for large $|s|$ either the root locus is close to a pole or it is asymptotic to the real axis.  

A plot with real $s$ shows that there are 2 values of $s$ for some values of $k$ (Figure \ref{fig1}).  This indicates a split  of the root locus on the negative real axis, with one branch going to $-\infty$ and the other to $0.$
\begin{figure}
\centering \includegraphics[width=2.5in]{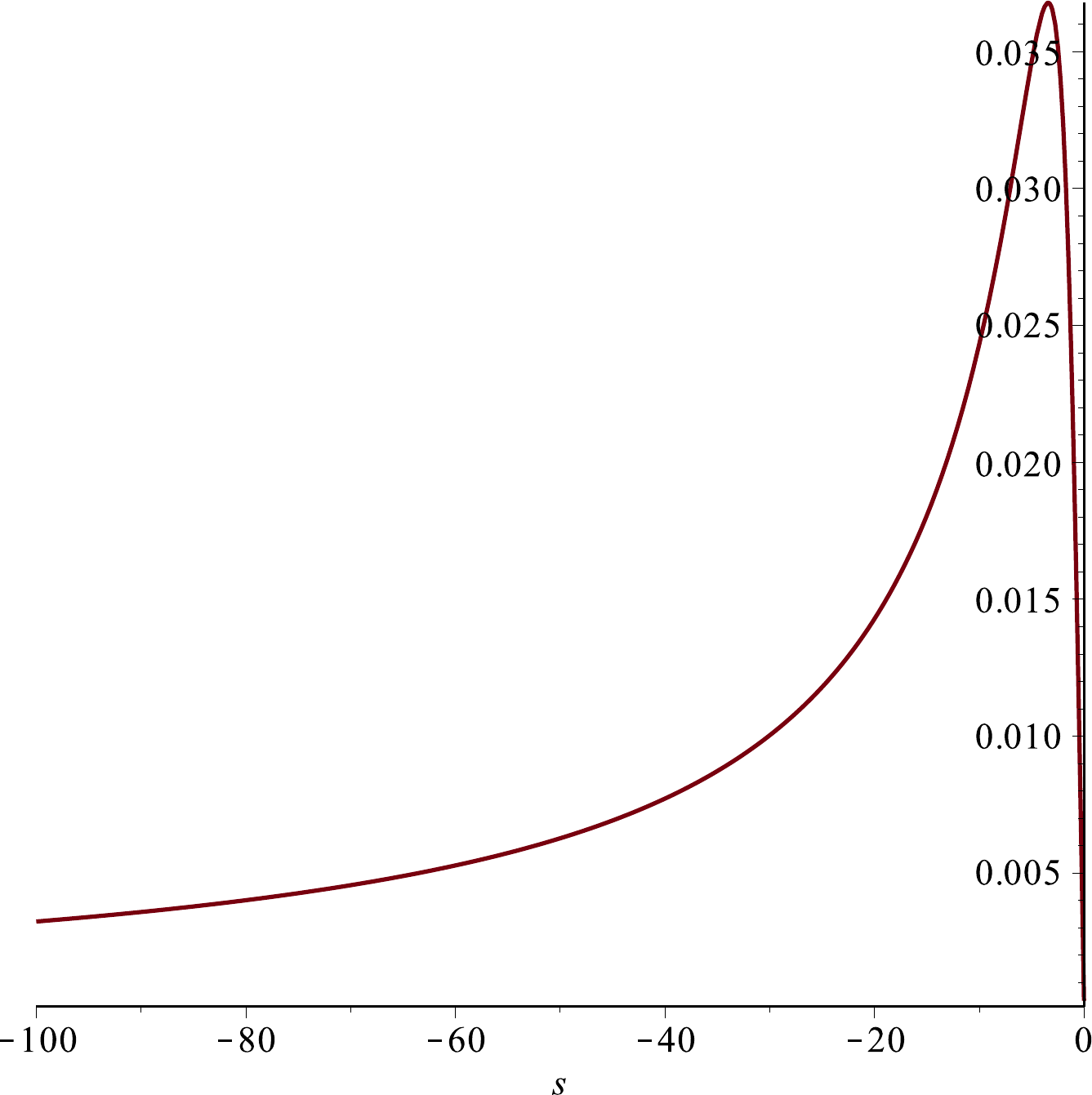}
\caption{$-G(s)$ for real $s$ for a simple wave equation. Since $-G(s)>0$ for all real $s$ and $\frac{1}{k} = - G(s)$ defines the root locus, the entire real axis is in the root locus, as proven in section 6. The plot  indicates that for some values of $k$ there are two branches of the root locus on the real axis. The value of $k$ at which the root locus reaches the real axis is determined by the maximum value of $G$. In this case, the root locus reaches the real axis when $k=27.2.$ }
\label{fig1}
\end{figure}
\begin{figure}
\centering \includegraphics[width=3.in,height=3.8in]{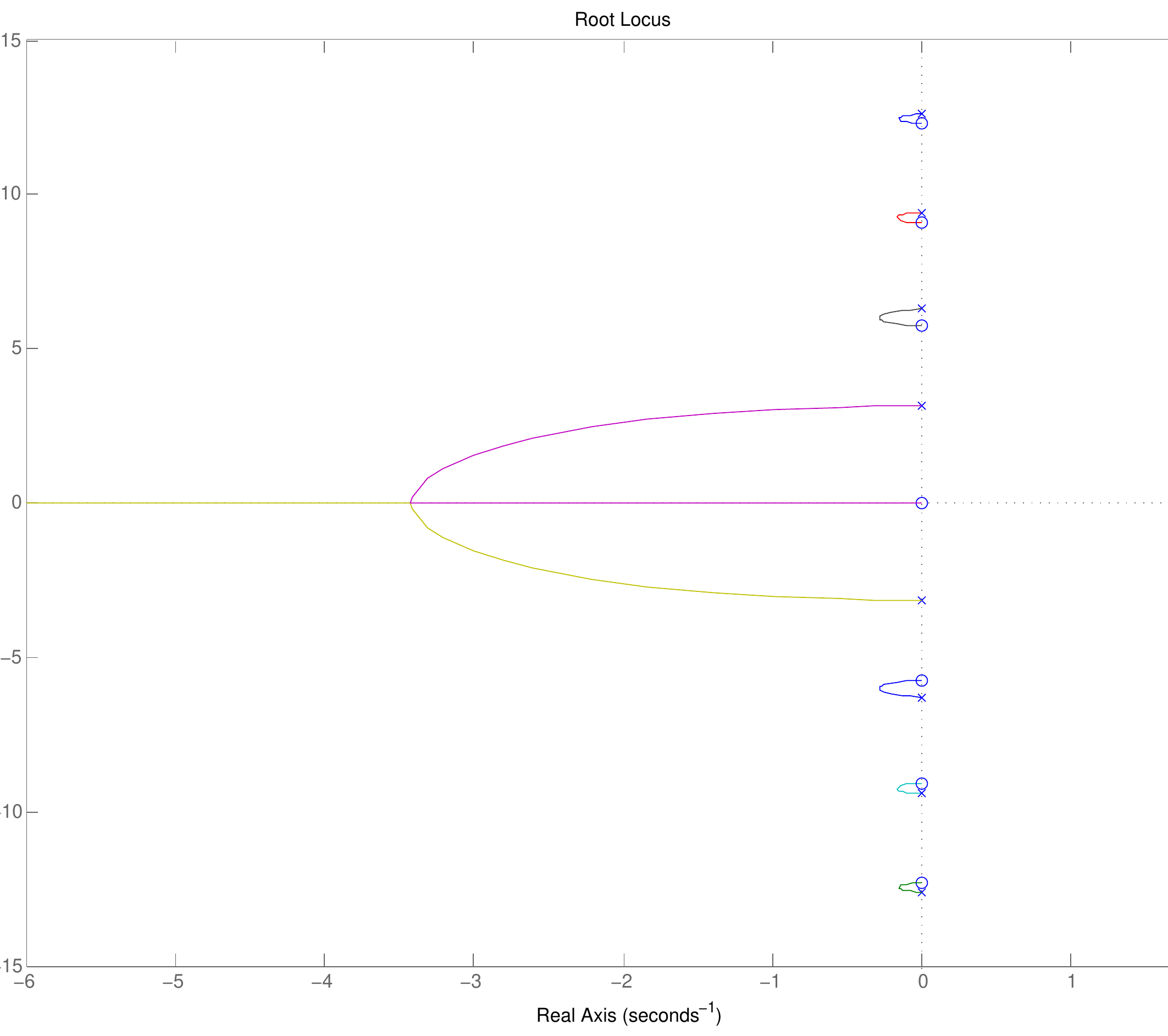}
\caption{Root locus of a system with the first 4 eigenvalues and 3 zeros of the wave equation example. Zeros are indicated with $\circ$ and eigenvalues by $x .$ Note that the zeros become very close to the  eigenvalues. The root locus starting at the smallest frequencies moves to the real axis and splits. The other branches move from the eigenvalues to the zeros}
\label{fig2}
\end{figure}
Thus,  the root locus is qualitatively quite similar to that of an analogous finite-dimensional system. There are  a number of branches curving from an eigenvalue to a zero, while two branches curve from eigenvalues to the real axis where they split. This is illustrated in Figure \ref{fig2}.
\end{eg}

\begin{eg} 
{\bf (Plate with Boundary Damping)}
Consider  vibrations in a plate or membrane  on  a bounded connected region $\Omega$ with
boundary $\Gamma.$
 The region $\Omega \subset \mathbb R^n $ has Lipschitz boundary $\Gamma$, where
  $\Omega$ is such that the embedding of $H^1 (\Omega)$ into $L^2 (\Omega)$
 is compact. Assume also $\Gamma =
\overline{ \Gamma_0 \cup \Gamma_1 }$
 and $\Gamma_0$, $\Gamma_1$ are disjoint open subsets of $\Gamma .$
  Consider the  control system description
\begin{equation}
\begin{array}{ll}
 \ddot{ w}    = \grad^2 w  + f u(t)  , \quad  & \Omega \times (0, \infty )  , \\
 w( x , 0)  = w_0, \; \dot{w} (x,0)   =  w_1, \quad & x \in  \Omega, \\
 w (x, t )  = 0, \quad  &  \Gamma_0 \times (0, \infty ), \\
 \frac{\partial w (x,t)}{\partial n} + \dot{w} (x,t)  = 0 &  \Gamma_1 \times (0, \infty ), \\
 y(t) = \int_{\Omega} f(x) \dot{w} (x,t) dx, & [0, \infty).
 \end{array}
 \label{wave}
 \end{equation}
Also assume that $f \in L^2 (\Omega)$ is  chosen so that  the system is approximately controllable/observable.
Define the self-adjoint operator $A_o$ on $L^2 (\Omega)$ by
\begin{align*}
A_o f &= -\grad^2 f, \\
 \dom(A_o) &= \{ f \in H^2 (\Omega) \cap H_{\Gamma_0}^1 (\Omega) |  \, \frac{ \partial f }{\partial n}|_{\Gamma_1} f = 0 \}.
\end{align*}
This leads to the abstract second-order differential equation 
$$ \ddot{z} (t)  + A_o z (t) = F u (t) $$
where the bounded operator $F$ is defined by $Fu=f(x) u.$
This is exactly in the class (\ref{2nd}) where $H =L^2 (\Omega )$  and $V=\dom(A_o^{\half} ) =H_{\Gamma_0}^1 ( \Omega) $ is the completion of $\dom(A_o)$ in
the norm $ (A_o z, z )^{1/2} .$  
It is well-known that this defines a contraction semigroup on $V \times H$; see \cite[e.g]{JMT} for details.

Although except for very special choices of $\Omega$ the transfer function cannot be calculated exactly, the results of this section ensure that, there is one zero at $0$, the zeros alternate with 
the eigenvalues on the imaginary axis and they asymptote towards the eigenvalues.  For all positive choices of gain $k$, the system is stable. The zero at $0$ is the terminus of a branch of the root locus and so there is a limit to the improvement in settling time that can be achieved with constant feedback. As for finite-dimensional systems, the system eventually becomes over-damped. 
\end{eg}

\section{Future Research}

In this paper, a rigorous definition of the root locus was provided for systems with bounded control and observation. Results on the qualitative nature of the root locus for collocated systems with a self-adjoint or skew-adjoint generator were also obtained.
Extension to systems with unbounded control and observation requires some care, since for unbounded feedback the spectrum can change dramatically, and a system with a complete set of eigenfunctions corresponding to an infinite sequence of eigenvalues can be perturbed to one with an empty spectrum. Since realistic models of actuators and sensors typically lead to bounded control and observation (see for instance \cite{JM12,ZLMV}) this is primarily of theoretical interest.

A significant open family of questions however is the qualitative nature of the root locus. The results in section 6 suggest that the root locus is in general similar to that of Example 6.9, shown in Figure 3, but this remains to be proven. Qualitative results for damped second-order systems, and for non-collocated systems are also desirable.

It is also shown in this paper that in many cases the invariant zeros are in the pseudo-spectrum of the eigenvalues. This has consequences for numerical calculations, in particular for order reduction.  Generalization of these results, if they do in fact generalize, points to the need for greater knowledge of the pseudo-spectrum of operators on infinite-dimensional spaces.

\bibliographystyle{IEEEtran}
\bibliography{rootlocus}

\end{document}